\documentclass[11pt,twoside,a4paper]{book}

\usepackage[T1]{fontenc}
\usepackage[brazil]{babel}
\usepackage[latin1]{inputenc}
\usepackage[pdftex]{graphicx}           
\usepackage{setspace}                   
\usepackage{indentfirst}                
\usepackage{makeidx}                    
\usepackage[nottoc]{tocbibind}          
\usepackage{courier}                    
\usepackage{type1cm}                    
\usepackage{listings}                   
\usepackage{titletoc}
\usepackage[fixlanguage]{babelbib}
\usepackage[font=small,format=plain,labelfont=bf,up,textfont=it,up]{caption}
\usepackage[usenames,svgnames,dvipsnames]{xcolor}
\usepackage[a4paper,top=2.54cm,bottom=2.0cm,left=2.0cm,right=2.54cm]{geometry} 
\usepackage[pdftex,plainpages=false,pdfpagelabels,pagebackref,colorlinks=true,citecolor=DarkGreen,linkcolor=NavyBlue,urlcolor=DarkRed,filecolor=green,bookmarksopen=true]{hyperref} 
\usepackage[all]{hypcap}                    
\fontsize{60}{62}\usefont{OT1}{cmr}{m}{n}{\selectfont}
\usepackage{amsmath}
\allowdisplaybreaks

\usepackage[final]{pdfpages} 
\fontsize{60}{62}\usefont{OT1}{cmr}{m}{n}{\selectfont}
\cleardoublepage
\normalsize

\begin{document}
\frontmatter 

\onehalfspacing  

\thispagestyle{empty}
\begin{center}
    \vspace*{2.3cm}
    \textbf{\Large{Aproximação do Equilíbrio e \\ Tempos Exponenciais para o \\ Passeio Aleatório no Hipercubo}}\\
    
    \vspace*{1.2cm}
    \Large{\textbf{Cláudia Peixoto}}
    
    \vskip 2cm
    \textsc{
    Dissertação apresentada\\[-0.25cm] 
    ao\\[-0.25cm]
    Instituto de Matemática e Estatística\\[-0.25cm]
    da\\[-0.25cm]
    Universidade de São Paulo\\[-0.25cm]
    para obtenção do grau de mestre\\[-0.25cm]
    em\\[-0.25cm]
    Estatística}
    
    \vskip 1.5cm
    Área de Concentração: \textbf{Probabilidade}\\
    Orientador: \textbf{Prof. Dr. Antonio Galves}

   	\vskip 1cm
    São Paulo, fevereiro de 1992
\end{center}
\newpage
\begin{center}
	\textbf{Abstract}
\end{center}
	We study a random walk in a N dimensional hypercube and exhibit results about stopping times when N diverges. The first theorem discusses the time in which two coupling processes spend to meet. A corollary provides a majorant for the velocity of convergence to equilibrium. Other three theorems treat, respectively, the time of first return to a point, the time of first return to a fixed set and the time of first arrival in a random set. We prove that these times, under a suitable rescaling, converge in law to a mean one exponential random time.

\mainmatter
\newpage

\def\eh{\rlap{e}\'{ }}
\def\ss{\smallskip\smallskip}
\baselineskip= 18 pt
\noindent{\bf I. Introdu\c c\~ao:}

\ss

\ss

\ss

\ss

Neste \ trabalho estudaremos \ tempos de parada para um passeio a\-le\-a\-t\'orio homog\^eneo em um hipercubo de dimens\~ao $ N $, obtendo resultados assint\'oticos.

\ss

\ss

A cada instante o processo assumir\'a uma configura\c c\~ao que ser\'a uma seq\"u\^encia de $ N $ elementos pertencentes ao conjunto $ \{ -1,+1 \} $ .
Teremos portanto $ 2^N $ configura\c c\~oes distintas que podem ser consideradas os v\'ertices do hipercubo. 

\ss

\ss

A evolu\c c\~ao do processo dar-se-\'a em tempo discreto e pode ser descrita da seguinte maneira: a cada passo com probabilidade $ { 1 \over 2} $ o passeio permanecer\'a no mesmo lugar e com probabilidade $ { 1 \over 2 } $ modificar\'a um de seus elementos escolhido de maneira uniforme.

\ss

\ss
Em nosso primeiro teorema, exibiremos a escala de tempo em que dois passeios aleat\'orios acoplados se encontram quando $ N $ diverge. O mesmo teorema pode ser encontrado em [2], mas aqui, a demonstra\c c\~ao est\'a bastante simplificada.

\ss

\ss

O segundo teorema trata do instante do primeiro retorno a uma posi\c c\~ao j\'a visitada pelo passeio. Neste teorema caracterizaremos, com probabilidade 1, como ser\'a este primeiro retorno quando $ N $ diverge. Este resultado faz parte de um artigo de \'Avila,Cassandro e Galves.

\ss

\ss   

O terceiro teorema trata do tempo de retorno a um conjunto fixado.  Este tempo, devidamente normalizado, converge em lei a uma exponencial de par\^ametro 1 quando $ N $ diverge. O que foi feito neste teorema generaliza um resultado de  Bellman e Harris  [3] , onde \eh\ tratado o tempo de retorno a uma \'unica posi\c c\~ao fixada. O artigo de Bellman e Harris estuda o modelo de Ehrenfest do qual o passeio aleat\'orio no hipercubo \eh\ uma esp\'ecie de vers\~ao microsc\'opica.  A nossa demonstra\c c\~ao aborda o problema de maneira an\'aloga a de Cassandro, Galves, Olivieri e Vares em [2] juntamente com o segundo teorema.  

\ss

Nosso quarto resultado refere-se ao instante de chegada do passeio  a um conjunto aleat\'orio $ \  M  \subset H_N \ $ . \  Cada ponto do hipercubo pertencer\'a a $ \ M \ $  com probabilidade $  { 1 \over N^{\gamma}}, \gamma>0 ,\ $
independentemente dos demais pontos e do passeio. Pontos em $ M $ ser\~ao chamados de pontos pretos.

\ss

\ss

Mostramos que o tempo necess\'ario para o passeio alcan\c car $ M $, normalizado pela densidade de pontos pretos, converge em lei a uma exponencial de par\^ametro 1 quando $ N $ diverge.

\ss

\ss

O modelo descrito acima foi estudado por Cassandro, Galves e Picco  em [2], o qual serviu de motiva\c c\~ao geral para este trabalho. Nosso objetivo foi desenvolver e refinar alguns resultados que ali aparecem.  

\ss

\ss

Na seq\"u\^encia apresentaremos a descri\c c\~ao do modelo e enunciaremos os quatro principais resultados. Ap\'os essa se\c c\~ao seguir\~ao as respectivas  demonstra\c c\~oes .

\includepdf[pages=1-last]{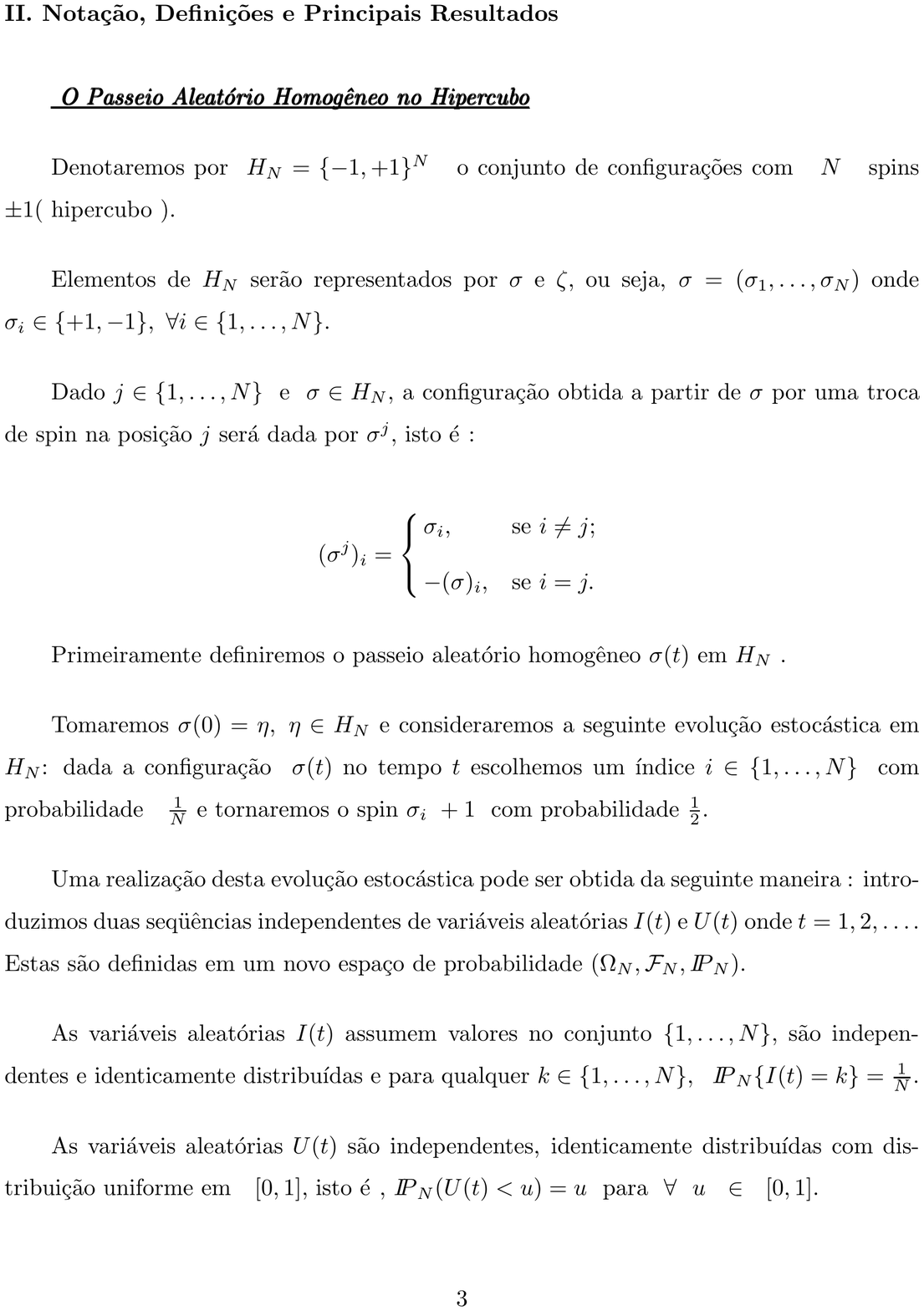} 

\def\eh{\rlap{e}\'{ }}
\def\Eh{\rlap{E}\'{ }}
\def\bbR{{I\kern-0.3emR}}
\def\bbP{{I\kern-0.3emP}}
\def\bbC{{I\kern-0.6emC}}
\def\bbN{{I\kern-0.3emN}}
\def\newend{\vbox to 7pt{\hbox to 7pt{\vrule height 7pt width 7pt}}}
\baselineskip = 18pt
\def\ss{\smallskip\smallskip}
\def\s{\smallskip}
\hyphenation {pro-ba-bi-li-da-de}
\hyphenation {a-le-a-to-ri-as}
\def\bff#1{\setbox0=\hbox{$#1$}
	\kern-.025em\copy0\kern-\wd0
	\kern.05em\copy0\kern-\wd0
	\kern-.025em\raise.0433em\box0 }

\noindent{\bf III. Resultados}

\ss

\ss

\ss

\ss

\noindent $ \bff { \hbox {III.1 \ \underbar{\it { Demonstração do Teorema I}}}} $

\ss

\ss

A seguir verificaremos que dois passeios homog\^eneos acoplados com configura\c c\~oes
iniciais tais que sua dist\^ancia seja m\'axima, encontrar-se-\~ao em
um tempo de ordem \ $ NlogN \ $ com probabilidade $ 1 \ $ quando
\ $ N $ \ diverge.

\ss

\ss

\ss

\ss

\ss

\noindent{\bf Teorema I: }

\ss

$$ \lim_{ N \rightarrow \infty} { t^{-}_N \over {N logN}} = 1  \hbox{ \ \    em 
	probabilidade.} $$

\ss

\ss

\noindent{\bf Demonstra\c c\~ao : }

\ss

\ss

Como mencionamos, $ \sigma^{+}(t) $ \ e \ $ \sigma^{-}(t) $ \ ser\~ao constru\'\i dos usando o mesmo  $\omega$ , isto \eh\  , a mesma escolha de \'\i ndices
$ I(t,\omega) $ e a mesma escolha de $ \ U(t,\omega) $ da seguinte maneira:

\ss

\ss
\ss

dados $ \sigma^{+}(t)  , \  \sigma^{-}(t) $ e $ I(t+1)=i, $

\ss

se $ U(t+1) < {1 \over 2} $ \ \ \ \ ent\~ao \  $ \sigma^{+}(t+1) = +1,  \ \sigma^{-}(t+1)= +1; $

\ss

\ss

se $ U(t+1) > { 1 \over 2} $ \ \ \ \ ent\~ao \ $ \sigma^{+}(t+1) = -1,  \ \sigma^{-}(t+1) = -1. $

\ss

\ss

\ss

\ss

Chamo \ $ D_N (n) \ = \ 
{ 1 \over {2N}} \ \sum\limits_{i=1}^N 
\ \vert \ \sigma^{+}_i (n) - \sigma^{-}_i (n) \  \vert 
$ \ , a dist\^ancia no instante $ n $  entre os dois passeios.

\ss

\s


Pela \ evolu\c c\~ao \  de $ \ \sigma(t) \  $  temos \ que \ $ D_N (n)  \ $ \ \eh\  \ uma  \ Cadeia \ de
\ Markov \ em $ \{ 0,{1 \over N}, \dots , {{N-1} \over N },1 \}, $ com probabilidade de transi\c c\~ao dada por:

\ss

\ss

$$ \bbP( x,y ) = \begin{cases}
1 - x , &se $ y=x$; \\
x , &se y=x-\frac{1}{N}; \\
0,& \text{caso contário.} \\
\end{cases} $$

\ss

Primeiramente notamos que :

$$ t^{-}_N = \sum\limits^{N}_{k=1} \lambda_k  \ , $$

\ss

\ss

onde $ \lambda_1 = 1 $ \ e \ $ \lambda_k $ \ para $ \ k = 2,\dots,N $ \ s\~ao
vari\'aveis aleat\'orias independentes, geom\'etricas de par\^ametro
$ \ p_k = 1 - { ( k - 1 ) \over N } . $

\ss

O tempo $ \ \lambda_k \ $ \eh\ exatamente o n\'umero de passos que o processo 
$ \ D_N(n) \ $ gasta para ir de $ 1 - { {k-1} \over N } \ $ a 
$ \ 1 - { k \over N } . $

\ss

\ss

\ss

Temos ent\~ao que suas esperan\c cas e vari\^ancias valem respectivamente :

\ss

\ss

$ E ( \lambda_k ) = { 1 \over p_k } = { N \over { N - k + 1 } },$

\ss

\ss

\ss

$ V ( \lambda_k ) = { q_k \over {p_k}^2 } = {{ N ( k - 1 )} \over { ( N - k + 1 )^2
	}} .$
	
\ss

\ss

Agora,
	
\begin{align*}
	\ \ \ \ \ \ \ \ \ \ \  E\bigl( t_{N}^{-} \bigr) &= \sum_{k=1}^{N} E\bigl( \lambda_k \bigr) \cr \cr
	&= \sum_{k=1}^{N} { N \over { N - k + 1 } } \cr \cr
	&= N \sum_{k=1}^{N} { 1 \over k } \cr
\end{align*}

Como   $ \ \ \  \int\limits^{N}_{1} {1 \over x } = logN $ , \ temos:

$$ N log \bigl(N-1\bigr) \leq E \bigl( t_N^{-} \bigr) \leq N \bigl( 1 + log N \bigr) .$$

\ss

\ss

\ss

Assim para qualquer $ \epsilon > 0 $,

\begin{align*}
 \bbP \Bigl( | t_N^{-} - E ( t_N^{-} ) | > 
 \epsilon  \ E ( t_N^{-} )
 \Bigr) &\leq { {V( t_N^{-} ) } \over { \epsilon^2 
 		[ E( t_N^{-} ) ]^2 }}
 \cr \cr
 &= { {\sum\limits^{N}_{k=1} { {N( k - 1)} \over { (N-k+1)^2 } } }
 	\over { \epsilon^2 [ \sum\limits^{N}_{k=1} {N \over { N-k+1 } }]^2 } }
 \cr \cr
 &= { { N \ \sum\limits^{N}_{k=1} { (k-1) \over {(N-k+1)^2 }} } \over
 	{ \epsilon^2 \ N^2 \bigl[ \ \sum\limits^{N}_{k=1} { 1 \over {(N-k+1)} } \bigr]^2
 	}} \cr \cr 
 	&\leq { { N(N-1) \bigl[ { 1 \over {(N-1)^2}} + \dots + 1 \bigr] } \over
	{ \epsilon^2 \ N^2 log^2 N }} \cr \cr 
 	&= { { N(N-1) \ \sum\limits^{N-1}_{k=1} 
			{ 1 \over {k^2} } } \over
 		{ \epsilon^2 \ N^2 log^2 N }} \cr \cr 
 	&\leq { { N ( N-1 )
 			\bigl[ 1 + \sum\limits^{N}_{k=2} { 1 \over {(k^2 - 1 ) }} \bigr]} \over
 		{ \epsilon^2 \ N^2 log^2 N } } \cr \cr 
 	&= { {N(N-1) \bigl[ 1 + {3 \over 4} - {{ (2N + 1)} \over {2N(N + 1)}} \bigr]
 		} \over { \epsilon^2 N^2 log^2 N }} .  \cr
\end{align*}
		
\ss
		

Passando ao limite quando \ $ N $ \ diverge temos que :

\ss

$$ \lim\limits_{N \rightarrow \infty} {\bbP \Bigl( 
			| {{t_N^{-}} \over
				{E(t^{-}_N) }} - 1 | > \epsilon \Bigr) }= 0 .$$

\ss

\hfill\newend

\ss
		
\ss

\noindent{\bf Corol\'ario I: } \ Sejam $ \ \sigma^{\eta} $ e $ \ \sigma^{\zeta} $ passeios homog\^eneos 
	em $ H_N $ constru\'\i dos simultaneamente, com configura\c c\~oes iniciais $ \eta, \zeta $. 
	Suponha que $ D_N(0) = { {[Nf]} \over N} $ , \ onde  $ \ 0 \leq f \leq 1 . $ 
	\ Ent\~ao,

		$$ \lim_{N \rightarrow \infty} \bbP\Bigl( \sigma^{\eta}( t(N) ) \not= \sigma^{\zeta}( t(N) ) \Bigr) = 0 , $$
		
		\ss

		para qualquer $ t(N) $ \ que satisfa\c ca $ \lim\limits_{N \rightarrow \infty} { {t(N)} \over { NlogN} } = \infty .$

		\ss
		
		\ss
		
		\noindent{\bf Demonstra\c c\~ao: } \ Primeiro notamos que:

		\ss
		
		$ t_N^{-} = \inf\bigl( t>0 : \{ I(1),\dots,I(t) \} = \{ 1,\dots,N \} \bigr) . $

		\ss
		
		Assim,
		
		\begin{align*}
		 \bbP\Bigl( \sigma^{\eta}( t(N) )  \not= \sigma^{\zeta}( t(N) ) \Bigr) &\leq \bbP( t_N^{-} > t(N) ) \cr \cr
		 &\leq { {N(logN + 1)} \over { t(N)}} . \cr
		\end{align*}
		
		\ss
		
		Portanto por hip\'otese temos que:
		
		$$ \lim_{N \rightarrow \infty} \bbP\Bigl( \sigma^{\eta}( t(N) ) \not=
		\sigma^{\zeta}( t(N) ) \Bigr) = 0 . $$

		\hfill\newend
		
		O pr\'oximo corol\'ario nos fornecer\'a um majorante para a velocidade de converg\^encia ao equil\'\i brio.
		
		\vfill\eject

		\noindent{\bf Corol\'ario II: } \ Seja $  \ t(N) $ tal que $ \lim\limits_{N \rightarrow \infty} 
		{ {t(N)} \over { NlogN} } = \infty .  $ \ \  Ent\~ao,
		
		\ss
		
		$$\lim_{N \rightarrow \infty} \Big\vert \bbP \Bigl( \sigma^{+}( t(N) ) = \zeta \Bigr) - { 1 \over {2^N} } \Big\vert = 0 , $$
		
		\ss
		
		onde $ \zeta \in  H_N $ .
		
		\ss

		\noindent{\bf Demonstra\c c\~ao: } 
		
		\ss
		
		Seja $ \ \nu(\eta) = { 1 \over {2^N} } $ a medida uniforme em $ H_N $. Esta medida \eh\ invariante com respeito \`a evolu\c c\~ao do passeio homog\^eneo, ou seja,

		$$ \nu(\zeta) = \sum_{\eta \in H_N} \nu(\eta) \bbP\Bigl( \sigma^{\eta}(t(N)) = \zeta \Bigr) .$$

		\ss

		Assim,
		
		\ss
		
		\begin{align*}
		 \big\vert \bbP\bigl(\sigma^{+}(t(N)) = \zeta\bigr) - { 1 \over {2^N} } \big\vert 
		 &= \big\vert \bbP\bigl(\sigma^{+}(t(N)) = \zeta\bigr) - \sum_{\eta \in H_N} \nu(\eta) \bbP\bigl(\sigma^{\eta}(t(N)) = \zeta\bigr) \big\vert \cr \cr
		 &\leq \sum_{\eta \in H_N} \nu(\eta) \big\vert \bbP\bigl(\sigma^{+}(t(N) ) = \zeta\bigr) - \bbP\bigl(\sigma^{\eta}(t(N))=\zeta\bigr) \big\vert  \cr \cr
		 &\leq \sum_{\eta \in H_N} \nu(\eta) \sup_{\eta} \bbP\bigl(\sigma^{+}(t(N)) \not= \sigma^{\eta}(t(N)) \bigr) \cr
		\end{align*}
		
		\ss
		
		\ss
		
		Passando ao limite em $ N $  e utilizando o  corol\'ario I a express\~ao acima vai a zero e portanto temos o resultado.

		\hfill\newend

		\ss
		
		\noindent{\bf Nota:} O processo $ \xi(t) $ \eh\ peri\'odico, al\'em disso, para qualquer acoplamento, dois processos com configura\c c\~oes iniciais que diferem em um n\'umero \'\i mpar de posi\c c\~oes nunca poder\~ao se encontrar, sendo que a dist\^ancia m\'\i nima entre eles ser\'a   a equivalente a uma posi\c c\~ao diferente.  Do teorema I, temos que a dist\^ancia em um tempo da ordem  $ NlogN $ entre dois passeios acoplados $ \xi^{\eta}(t) $ e $ \xi^{\zeta}(t) $ ser\'a menor ou igual a $ { 1 \over N } $ quando $ N $ diverge. 
		
		\ss
		
		\newpage
		
\def\eh{\rlap{e}\'{ }}
\def\Eh{\rlap{E}\'{ }}
\def\bbR{{I\kern-0.3emR}}
\def\bbP{{I\kern-0.3emP}}
\def\bbC{{I\kern-0.6emC}}
\def\bbN{{I\kern-0.3emN}}
\def\newend{\vbox to 7pt{\hbox to 7pt{\vrule height 7pt width 7pt}}}
\baselineskip = 18pt
\def\ss{\smallskip\smallskip}
\def\s{\smallskip}
\hyphenation {pro-ba-bi-li-da-de}
\hyphenation {a-le-a-to-ri-as}
\def\bff#1{\setbox0=\hbox{$#1$}
	\kern-.025em\copy0\kern-\wd0
	\kern.05em\copy0\kern-\wd0
	\kern-.025em\raise.0433em\box0 }

$\bff{ \hbox{ \noindent{\bf III.2 }\underbar{\sl { Demonstra\c c\~ao do Teorema II  }}}} $

\ss

\ss

Considere os seguintes eventos:

\ss

\ss

$ \ \  J_1 = \Bigl\{ (i_1,i_2) \in \{1,\dots,N\}^{2} : i_1 = i_2 \Bigr\} $

\ss

\ss

\ss

$ \ \ J_l = \Bigl\{ (i_1,\dots,i_{2l} ) \in \{1,\dots,N\}^{2l} :
\sum\limits_{k=1}^{2l} {\bf 1}_{ \{ i=i_k\} } \in \{ 0,2,\dots,2l\} \ \forall i \in \{1,\dots,N\} $

\ \ \ \ \ \ \ \ \ \ \  e
$ \ \bigl( i_{k},\dots,i_{k+2m-1} \bigr) \not\in J_m $, $ \forall m  \in \{1,\dots,l-1\}, \  k \geq 1 \hbox{ \ e \ } k+2m-1 \leq 2l . $

\ss

\ss

\ss

\noindent{\bf Lema I: } \ \  Para $ l \geq 3 $ temos:

\ss
$$ \bbP\Bigl( \{ I(1),\dots,I(2l) \} \in J_l \Bigr) \leq { 8 \over {N^3}} . $$

\ss

\noindent{\bf Demontra\c c\~ao: } 

$$\bbP \Bigl( \{ I(1),\dots,I(2l) \} \in J_l \Bigr) = { {\big\vert J_l \big\vert} \over { N^{2l}} } . $$

\ss

Primeiramente notamos que $ \big\vert J_l \big\vert \leq 2 N^{2l-2} $ pois nas duas \'ultimas posi\c c\~oes os \'\i ndices est\~ao fixados a menos de uma permuta\c c\~ao.

\ss

\ss

Agora dividiremos o conjunto $ J_l $ em dois conjuntos disjuntos: aqueles em que as tr\^es \'ultimas posi\c c\~oes s\~ao ocupadas por \'\i ndices distintos entre si e aqueles em que nas tr\^es \'ultimas posi\c c\~oes aparecem apenas dois \'\i ndices distintos entre si.

\ss

\ss

Denotaremos esses conjuntos por $ J^{\prime}_l $ e $ J^{\prime\prime}_l $ respectivamente.

\ss

\ss

Temos ent\~ao que: $ \ \ \big\vert J_l \big\vert = \big\vert J^{\prime}_l \big\vert + \big\vert J^{\prime\prime}_l \big\vert . $ 

\ss


Note que:

\ss

$ \ \ \hbox{a) \ } \  \big\vert J^{\prime}_l \big\vert \leq 3! \ N^{2l-3} $ , pois as tr\^es \'ultimas posi\c c\~oes devem estar fixadas, a menos de uma permuta\c c\~ao, para que ocorra retorno.

\ss

\ss

$ \ \  \hbox{ b) \ } \ \big\vert J^{\prime\prime}_l \big\vert \ \leq N  \ \big\vert J_{l-1} \big\vert $ , pois eliminando o par que aparece nas tr\^es \'ultimas posi\c c\~oes temos exatamente um retorno do tipo $ J_{l-1} $; al\'em disso, o algarismo repetido pode assumir no m\'aximo $ N $ valores.

\ss

\ss

Portanto,

\begin{align*}
{{\big\vert {J}_l \big\vert} \over { N^{2l}}} &\leq
{ {3! \ N^{2l-3}} \over { N^{2l}} } + { {N \  |J_{l-1}| } \over { N^{2l}} } \cr \cr
&\leq { {3! \ N^{2l-3}} \over { N^{2l}}} + { {N \  2N^{2l-4} } \over {N^{2l}} }
\cr \cr
&\leq { 8 \over {N^3}} . \cr
\end{align*}

\ss

\hfill\newend

Introduziremos agora as seguintes vari\'aveis aleat\'orias:

\ss

\ss

$  \ \ S_N = \inf( k \geq 2 : \xi(k) \in V[0,k-1] ), \hbox{onde} \ V[0,k-1]=\{\xi(0),\dots,\xi(k-1)\};  $

\ss

$ \ \ \ \Gamma_1 = \inf ( k \geq 2 : I(k) = I(k-1) ),$

\ss

$ \ \ \ \Gamma_l = \inf ( k \geq 2l : ( I(k-2l+1), \dots, I(k) ) \in J_l ) . $

\ss

\ss
\ss

\ss
\noindent{\bf Lema II : } \ Para $ l \geq 3 $ temos:

$$ \bbP ( \Gamma_l \leq n ) \ \leq \ n  \ {  8 \over {N^3} } .$$

\noindent{\bf Demonstra\c c\~ao:}


\begin{align*}
\bbP(\Gamma_l \leq n) &= \sum_{k=2l}^{n} \bbP(\Gamma_l = k) \cr \cr
&= \sum_{k=2l}^{n} \bbP\Bigl( (I(j-2l+1),\dots,I(j)) \not\in J_l, j<k; \cr 
&\ \ \ \ \ \ \ \ \ \ \ \  (I(k-2l+1),\dots,I(k)) \in J_l \Bigr)  \cr \cr
&\leq \sum_{k=2l}^{n} \bbP \Bigl( (I(k-2l+1),\dots,I(k)) \in J_l \Bigr)  \cr \cr
&= (n-2 l) \bbP\Bigl( (I(1),\dots,I(2l)) \in J_l \Bigr)  \cr
\end{align*}

\ss

Portanto, utilizando o resultado do lema I temos:

$$ \bbP \Bigl( \Gamma_l \leq n \Bigr) \leq {{ 8 n } \over { N^3} } . $$

\hfill\newend

\ss
\ss

\noindent{\bf Teorema II: } \ Para $ S_N $ e $\Gamma_1 $ como definidos anteriormente temos,

\ss

$$ \lim_{N \rightarrow \infty} \bbP\Bigl( S_N = \Gamma_1 \Bigr) = 1. $$

\ss

\noindent{\bf Demonstra\c c\~ao: } \ Inicialmente observamos que $ S_N = 
\min\limits_{l \geq 1} \Gamma_l . $

\ss

\ss

Para provarmos este teorema \eh\ suficiente mostrar que:

\s

$$\lim_{N \rightarrow \infty} \bbP( \Gamma_1 < N^{1 + \delta} < \min_{{{N^{1+\delta}} \over 2} \geq l \geq 2} \Gamma_l )=1 , \hbox{ \ para algum }   \ \delta>0. $$


Primeiro notamos que :

\ss

\begin{align*}
\lim_{N \rightarrow \infty} \bbP(\Gamma_1 < N^{1 + \delta}) &= 
\lim_{N \rightarrow \infty} 1 - \bbP( \Gamma_1 > N^{1 + \delta} ) \cr \cr
&= \lim_{N \rightarrow \infty}  1 - ( 1 - { 1 \over N })^{N^{1+\delta}}  = 1 \cr
\end{align*}

\ss

Por outro lado,

$$\bbP( N^{1 + \delta} < \min_{{ {N^{1 + \delta} } \over 2 } \geq l \geq 2} \Gamma_l ) \geq 1 - { {2 N^{1+\delta}} \over {N^2}}
- \sum_{l \geq 3}^{{N^{1 + \delta} \over 2 }} { {8 N^{1 + \delta}} \over {N^{3}}}. $$

\ss

Se tomarmos $ 0<\delta<{ 1 \over 2} $ a \'ultima express\~ao  vai a 1 quando $ N $ diverge. Isto conclui a demonstra\c c\~ao do teorema.

\ss

\hfill\newend

\ss

\ss

\noindent{\bf Corol\'ario III:} Para $ 0< \gamma <1 $ temos: 

$$ \lim_{N \rightarrow \infty} \bbP\Bigl( | \{ \xi(0),\dots,\xi(N^{\gamma}t) \} | = N^{\gamma}t+1 \Bigr) = 1 . $$

\ss

\noindent{\bf Demonstra\c c\~ao:}

\ss

\begin{align*}
 \lim_{N \rightarrow \infty} \bbP\Bigl( | \{ \xi(0),\dots,\xi(N^{\gamma}t) \} | = N^{\gamma}t+1 \Bigr) &=
 \lim_{N \rightarrow \infty} \bbP\Bigl( S_N > N^{\gamma}t \Bigr) \cr \cr
 &= \lim_{N \rightarrow \infty} \bbP\Bigl( \Gamma_1 > N^{\gamma}t \Bigr) \cr \cr
 &= \lim_{N \rightarrow \infty} \Bigl( 1 - {1 \over N } \Bigr)^{N^{\gamma}t} = 1. \cr
\end{align*}

\ss

\hfill\newend

\ss

\vfill\eject

\noindent{\bf Corol\'ario IV: } \ A vari\'avel aleat\'oria $ \ N^{-1} S_N  $ \ converge em lei para a distribui\c c\~ao
exponencial de par\^ametro 1 .

\ss

\ss

\noindent{\bf Demonstra\c c\~ao:} \ Basta provarmos que para qualquer $ t>0, $ temos:

$$ \lim\limits_{N \rightarrow \infty} \bbP( S_N > t N ) = e^{-t} . $$

\s

Agora,

$$ \big\vert \bbP\Bigl( S_N > t N \Bigr) - \bbP\Bigl( \Gamma_1 > t N \Bigr) \big\vert \leq \bbP\Bigl( \Gamma_1 \not= S_N \Bigr). $$

\ss

\ss

Pelo teorema III, $ \ \lim\limits_{N \rightarrow \infty} \bbP\Bigl( \Gamma_1 \not= S_N \Bigr) = 0 . $

\ss

\ss

Por outro lado,

\ss

\ss

$$ \lim_{N \rightarrow \infty} \bbP\Bigl( \Gamma_1 > t N \Bigr) = \lim_{N \rightarrow \infty} \Bigl( 1 - { 1 \over N } \Bigr)^{Nt} = e^{-t}. $$

\ss

Portanto,

$$ \lim_{N \rightarrow \infty} \bbP\Bigl( S_N > t N \Bigr) = e^{-t} . $$

\hfill\newend

\newpage

\def\eh{\rlap{e}\'{ }}
\def\Eh{\rlap{E}\'{ }}
\def\bbR{{I\kern-0.3emR}}
\def\bbP{{I\kern-0.3emP}}
\def\bbC{{I\kern-0.6emC}}
\def\bbN{{I\kern-0.3emN}}
\def\newend{\vbox to 7pt{\hbox to 7pt{\vrule height 7pt width 7pt}}}
\baselineskip = 18pt
\def\ss{\smallskip\smallskip}
\def\s{\smallskip}
\hyphenation {pro-ba-bi-li-da-de}
\hyphenation {a-le-a-to-ri-as}
\def\bff#1{\setbox0=\hbox{$#1$}
	\kern-.025em\copy0\kern-\wd0
	\kern.05em\copy0\kern-\wd0
	\kern-.025em\raise.0433em\box0 }

\ss

\noindent $ \bff{ \hbox{\underbar{\sl {III.3 Demonstra\c c\~ao do Teorema III}}}} $

\ss

Obteremos agora resultados sobre o tempo de retorno do processo $ \sigma^{+}(t) $ ao conjunto $ V[0,N^{\gamma}] $.  A partir de agora $ V[0,N^{\gamma}] $ ser\'a denotado por $ F $ .   
\ss

Considere

$$ R_N = \inf \Bigl( t > N^{\gamma} : \sigma^{+}(t) \in F \Bigr), \hbox{ \ e } $$
$$ R^{\eta}_N = \inf \Bigl(t>0 : \sigma^{\eta}(t) \in F \Bigr). $$

\s

\ss

\ss

\noindent{\bf Proposi\c c\~ao I : } \  \ Para \ $ 0 < \gamma < 1 $ \ e \
$ 0 < \delta < { 1 \over 2} , $ temos:

$$ \lim_{N \rightarrow \infty} \bbP\Bigl(R_N > N^{1 + \delta} \Bigr) = 1 . $$

\noindent{\bf Demonstra\c c\~ao: }

\ss

\ss

Do teorema II temos que:

$$ \lim_{N \rightarrow \infty} \bbP \bigl( \cup_{l \geq 2} \Gamma_l < N^{1+\delta}  \bigr) = 0. $$

Sendo assim, 

\begin{align*}
 \lim_{ N \rightarrow \infty} \bbP \bigl( R_N \leq N^{1+\delta} \bigr) &= \lim_{N \rightarrow \infty} \bbP \bigl( \sigma^{+}(N^{\gamma}+1) =
 \sigma^{+}(N^{\gamma}) \bigr) \cr \cr
 &\leq \lim_{N \rightarrow \infty}{ {N^{\gamma}} \over N } .\cr
\end{align*}

\ss

Como por hip\'otese $ 0<\gamma<1 $ o limite acima \eh\ zero.

\ss

Portanto,

$$\lim_{N \rightarrow \infty} \bbP\Bigl( R_N > N^{1 + \delta} \Bigr) = 1 .
$$

\hfill\newend

\ss

\vfill\eject

\noindent{\bf Proposi\c c\~ao II :} \ Seja $ F $ um conjunto fixado de cardinal $ N^{\gamma} , 0 < \gamma < 1 . $ Para $ \delta $ tal que

$ \ \ \ \ \ \ \ \ \ \ \ \  \  \ \gamma + \delta < 1 , \ \forall \ \eta \not\in F $ fixado temos:

$$ \lim_{N \rightarrow \infty} \bbP\bigl( R^{\eta}_{N} > N^{1+\delta} \bigr) = 1 . $$

\noindent{\bf Demonstra\c c\~ao:} \ Sem perda de generalidade a demonstra\c c\~ao ser\'a feita para o processo $ \xi^{\eta} $ . \ Primeiramente mostraremos que para uma configura\c c\~ao fixada em $ F $, o processo $ \xi^{\eta} $ gastar\'a um tempo maior que $ N^{1+\delta} $ para encontr\'a-la quando $ N $ diverge.

\ss

\ss

Fixe $ \zeta \in F $ . \  Defina  $ d(n) = { 1 \over 2}
\sum\limits_{i=1}^{N} | \ \xi^{\eta}_i(n) - \zeta_i  \ | $  \  a   dist\^ancia do passeio  \ $ \xi^{\eta} $ \ a $ \ \zeta \ $ no instante $ n $. \ Por hip\'otese, $ d(0) \geq 1 $ pois $ \zeta \not\in F $ .

\ss

Observe que $ d(n) $ evolue como o modelo de Ehrenfest; ou seja, uma cadeia de Markov com espa\c co de estados $ \{ 0,1,\dots,N \} $ \ e probabilidades de transi\c c\~ao dadas por:

$$  \bbP_{i,j} = \begin{cases}
 {{N-i}\over N}, &se  j=i+1;  \cr \cr
 { i \over N },  &se   j=i-1. 
\end{cases}$$ 

\ss

Sejam   

$ \ \ \ \ \  n_1 = \inf( n>1 : d(n)=2 ), $

$ \ \ \ \ \  n_2 = \inf(n>n_1 : d(n)=2 ) , $ e assim sucessivamente.

\ss

A cada retorno ao ponto "2", a probabilidade de em seguida visitar o ponto "0" antes de voltar ao "2" \eh\ $ {2 \over {N^2}} $.

\ss

Seja $ K = \inf( k\geq1 : d(n_k + 2)=0 ), $ ent\~ao o tempo para $ \xi^{\eta} $ alcan\c car $ \zeta $ \eh\  maior ou igual a $ K $.

\ss

Mas,
\begin{align*}
\bbP\bigl( K>t \bigr) &= \sum_{j=1}^{\infty} \bbP\bigl(K=t+j\bigr) \cr \cr
&= \sum_{j=1}^{\infty} \bbP\bigl( d(n_1+2) \not= 0, \dots, d(n_{t+j-1}+2)\not= 0, d(n_{t+j}+2)=0 \bigr) \cr \cr
&= \sum_{j=1}^{\infty} \bigl( 1 - {2 \over {N^2}} \bigr)^{t+j-1} { 2 \over {N^2}} \cr \cr
&= \bigl( 1 - {2 \over {N^2}} \bigr)^{t} . \cr
\end{align*}

Note que para $ t=N^{1+\delta} $ o limite acima vai a um quando $ N $ diverge.

\ss

Agora terminaremos a demonstra\c c\~ao observando que:
\begin{align*}
\bbP\bigl( \xi^{\eta}(u) &\in F, \hbox{ \ para algum \ } u \leq N^{1+\delta} \bigr) = \cr \cr
&\leq \sum_{\zeta \in F} \bbP\bigl( d(u)=0, \hbox{ para algum } \ u \leq N^{1+\delta} \ | \ d(0)= d(\zeta,\eta) \bigr) \cr \cr
&\leq N^{\gamma} \bbP\bigl( d(u)=0, \hbox{ para algum} \ u \leq N^{1+\delta} \ | \ d(0)=1 \bigr)  \cr \cr
&\leq N^{\gamma} \Bigl( 1 - \bigl( 1 - {2 \over {N^2}}\bigr)^{N^{1+\delta}} \Bigr) \cr \cr
&= N^{\gamma} \Bigl( 1 - exp\{ N^{1+\delta} log( 1 - {2 \over {N^2}}) \} \Bigr) \cr \cr
&= N^{\gamma} \Bigl( 1 - exp \{ N^{1+\delta} [ - ( { 2 \over {N^2}} +
{4 \over {2N^4}} + \dots ) ] \} \Bigr) \cr 
&= N^{\gamma} \Bigl( 1 - exp \{ -2N^{\delta-1} - 2N^{\delta-3} - {8 \over 3} N^{\delta - 5} - \dots \} \Bigr) \cr \cr
&=N^{\gamma} \Bigl( 1 - [ 1 - 2N^{\delta-1} + 2N^{2 \delta - 2} -{ 8 \over 6} N^{3 \delta - 3} + 4N^{4 \delta - 4} - \dots ] \Bigr) + o(N) \cr \cr
&= N^{\gamma} 2 N^{\delta-1} + o(N) . \cr
\end{align*}

Portanto passando ao limite e utilizando a hip\'otese $ \gamma+\delta<1 $ temos o resultado.

\hfill\newend

\noindent{\bf Proposi\c c\~ao III : } \ Para  $ 0<\gamma<1 $ \ e qualquer $ t(N) $ tal que 
$ \lim\limits_{N \rightarrow \infty} { {t(N) N^{\gamma} } \over { 2^N} } = 0 $ \  temos :

$$ \lim_{N \rightarrow \infty} \bbP \Bigl(
\ R_N > t(N) \ \Bigr) \ = \ 1 .$$

\noindent{\bf Demonstra\c c\~ao : }
\begin{align*}
\bbP\Bigl( R_N &\leq t(N) \Bigr) = \bbP\Bigl( \sigma^{+}(t) \in F,
\ \hbox{para algum} \ t\leq t(N) \Bigr) \cr \cr
&\leq \bbP\Bigl( \sigma^{+}(t) \in F,
\ \hbox{para algum} \ t \leq N^{1 + \delta} \Bigr) + \cr \cr
&+ \bbP\Bigl( \sigma^{+}(t) \in F, \ \hbox{para algum} \ N^{1 + \delta} < t \leq t(N) \Bigr)
\cr \cr
&= \bbP\Bigl( R_N < N^{1 + \delta} \Bigr) + \bbP\Bigl( \sigma^{+}(t) \in F, \ \hbox{para algum} \ N^{1+\delta}<t\leq t(N) \Bigr). \cr
\end{align*}

Como pela proposi\c c\~ao I $ \lim\limits_{N \rightarrow \infty} \bbP\Bigl(  R_N \leq N^{1+\delta} \Bigr) = 0 , $ temos:

\begin{align*}
\bbP\Bigl(R_N &\leq t(N) \Bigr) \leq o(N) + \bbP\Bigl( \sigma^{+}(t) \in F, \ \hbox{para algum} \ N^{1+\delta} < t \leq t(N) \Bigr) \cr \cr
&= o(N) + \bbP\bigl(\sigma^{+}(t) \in F, \hbox{ \ para algum \ } N^{1+\delta}<t<t(N) \bigr) \ -
\cr \cr
&- \bbP\bigl( \sigma^{\eta}(t) \in F, \hbox{\ para algum \ } N^{1+\delta}<t \leq t(N) \bigr) \ + \cr \cr
&+ \bbP\bigl( \sigma^{\eta}(t) \in F, \hbox{ \ para algum \ } N^{1+\delta}< t \leq t(N)
\bigr) \cr \cr
&\leq  o(N) + \sum_{\eta \in H_N} \nu(\eta) \Big\vert \bbP\bigl( \sigma^{+}(t) \in F, \hbox{ para algum \ } N^{1+\delta}<t\leq t(N) \bigr) \ - \cr \cr
&-\bbP\bigl( \sigma^{\eta}(t) \in F, \hbox{para algum} \ N^{1+\delta}<t \leq t(N) \bigr) \Big\vert \ + \cr \cr
&+ \sum_{\eta \in H_N} \nu(\eta) \bbP\bigl( \sigma^{\eta}(t) \in F, \hbox{ para algum \ } N^{1+\delta}<t \leq t(N) \bigr) \cr \cr
&\leq  o(N) + \sum_{\eta \in H_N} \nu(\eta) \sup_{\eta} \bbP\bigl(\sigma^{+}(N^{1+\delta}) \not= \sigma^{\eta}(N^{1+\delta}) \bigr) \ + \cr \cr
&+ \sum_{\eta \in H_N} \nu(\eta) \sum_{u=N^{1+\delta}}^{t(N)} \bbP\bigl( \sigma^{\eta}(t) \in F \bigr) \cr \cr
&\leq o(N) + \sum_{\eta \in H_N} \nu(\eta)  \sup_{\eta} \bbP\bigl( \sigma^{+}(N^{1+\delta}) \not= \sigma^{\eta}(N^{1+\delta}) \bigr) \ + \cr \cr
&+ { {(N^{\gamma}+1) t(N)} \over {2^N} } . \cr
\end{align*}

Passando ao limite em $ N $ , utilizando o teorema I e a hip\'otese, temos que a probabilidade acima vai a zero para algum $ 0<\delta<1. $

\ss


Portanto, 

$$\lim_{N \rightarrow \infty} \bbP\Bigl( R_N > t(N) \Bigr) = 1. $$

\hfill\newend

\ss

\noindent{\bf Proposi\c c\~ao IV :} Para  $\ \delta, \epsilon>0 $ e qualquer $ t(N) $ tal que 
\ $ \lim\limits_{N \rightarrow \infty} { { {t(N)}^{1-\epsilon} \nu(F) } \over { N(logN +1)} } = \infty $ \ temos :
$$ \lim_{N \rightarrow \infty} \bbP\bigl( R_N \leq t(N) \bigr) = 1 . $$

\s
\noindent{\bf Demonstra\c c\~ao:} \ Considere os eventos $ A_s = \{ \sigma(s) \in F \} $ e defina $ Z $ uma vari\'avel aleat\'oria positiva da seguinte maneira: 
$ Z=\sum\limits_{s=0}^{t(N)} {\bf 1}_{ \{A_s\} } .$ 
\ss

\ss

Note que  $ \{ Z>0 \} = \{ R^{\eta}_N \leq t(N) \} $ .

\ss
\ss

Aplicando desigualdade de Cauchy-Schwarz  ao produto \ $ Z {\bf 1}_{ \{ Z>0 \} } $ \ temos que:

$$ \bbP\bigl( Z>0 \bigr) \geq { {[E(Z)]^2} \over { E(Z^2)} } . $$
Assim,
\begin{align*}
&\bbP\bigl( Z > 0 \bigr) = \bbP\bigr( R^{\eta}_N \leq t(N) \bigr)  \cr \cr
&\geq { {[ E( \sum_{s=N^{1+\delta}}^{t(N)} {\bf 1}_{ \{A_s\} } ) ]^2 } \over
	{ [ E(\sum_{s=N^{1+\delta}}^{t(N)} {\bf 1}_{ \{A_s\} } )^2 ] } } + o(N) \cr \cr
&={{ [ \sum_{s=N^{1+\delta}}^{t(N)} \bbP(A_s) ]^2 } \over
	{ E( \sum_{s=N^{1+\delta}}^{t(N)} {\bf 1 }^2_{ \{A_s\}} + \sum_{u\not= s} {\bf 1}_{ \{ A_u\} } { \bf 1}_{ \{ A_s \} } )} }+ o(N) \cr \cr
&= { { (t(N)-N^{1+\delta} )^2 \nu(F)^2  } \over
	{ (t(N)-N^{1+\delta}) \nu(F) + \sum_{u \not= s} \bbP(A_u \cap A_s ) } } + o(N) \cr \cr
&= { {(t(N)-N^{1+\delta} )^2 \nu(F)^2 } \over
	{ (t(N)-N^{1+\delta}) \nu(F) + \sum_{|u-s|>N^{1+\delta}} \bbP(A_u \cap A_s)
		+ \sum_{|u-s| \leq N^{1+\delta}} \bbP(A_u \cap A_s ) } } + o(N)
. \cr 
\end{align*}

\vfill\eject
Agora observe que:

\noindent{\bf a)}
\begin{align*}
\sum\limits_{|u-s|>N^{1+\delta}}  \bbP\bigl( A_u \cap A_s \bigr)&=
\sum\limits_{k=N^{1+\delta}}^{t(N)} (t(N)-k+1) \bbP\bigl( \sigma(k) \in F, \sigma(0) \in F \bigr) \cr \cr
&= \sum_{k=N^{1+\delta}}^{t(N)} (t(N)-k+1) \sum_{\eta \in F} \nu(\eta) \bbP\bigl( \sigma(k) \in F | \sigma(0)=\eta \bigr) \cr \cr
&= \sum_{k=N^{1+\delta}}^{t(N)} (t(N)-k+1) \times \cr \cr
&\times \Bigl\{ \sum_{\eta \in F} \nu(\eta) \sum_{\xi \in H_N} \nu(\xi)  \Bigl[ \bbP\bigl( \sigma^{\eta}(k) \in F | \sigma(0)=\eta \bigr) - \cr \cr
&-\bbP\bigl( \sigma^{\xi}(k) \in F | \sigma(0)=\xi \bigr) \Bigr]  + \cr \cr
&+ \sum_{\eta \in F} \nu(\eta) \sum_{\xi \in H_N} \nu(\xi) \bbP\bigl( \sigma(k) \in F | \sigma(0)=\xi \bigr) \Bigr\}\cr \cr
&\leq \sum_{k=N^{1+\delta}}^{t(N)} (t(N)-k+1) \times
\cr \cr
&\times \Bigl[ \sum_{\eta \in F} \nu(\eta) \sum_{\xi \in H_N} \nu(\xi) \bbP\bigl( \sigma^{\eta}(k) \not= \sigma^{\xi}(k) \bigr) + \nu(F)^{2} \Bigr] \cr
\end{align*}

\ss
\vfill\eject
Usando a desigualdade de Markov e a majora\c c\~ao que aparece no teorema I obtemos que:

\begin{align*}
\sum_{|u-s|>N^{1+\delta}}\bbP\bigl(A_u \cap A_s \bigr) &\leq 
\sum_{k=N^{1+\delta}}^{t(N)} \bigl( t(N)-k+1 \bigr)  \times \cr \cr
&\times \bigl[ \sum_{\eta \in F} \nu(\eta) \sum_{\xi \in H_N} \nu(\xi) { {N(logN +1)} \over k }+\nu(F)^2 \bigr] \cr \cr
&=
\sum_{k=N^{1+\delta}}^{t(N)} \bigl( t(N) -k+1\bigr) \bigl[ \nu(F) { {N(logN+1)} \over k} + \nu(F)^2 \bigr] \cr \cr
&\leq \ t(N)log t(N) \ N(logN +1)\ \nu(F) + t(N)^{2} \nu(F)^{2}  \cr \cr
&\leq t(N)^{1+\epsilon} N(logN +1) \nu(F) + t(N)^2\nu(F)^2  . \cr
\end{align*}

\ss

\ss

\ss

\noindent {\bf b) } $$  \sum_{|u-s|<N^{1+\delta}}  \bbP\bigl( A_u \cap A_s \bigr) \leq t(N) N^{1+\delta} \nu(F) . $$

\ss

\ss

Passando ao limite em $ N $ temos por hip\'otese que:

$$ \lim_{N \rightarrow \infty} \bbP\bigl( R^{\eta}_N \leq t(N) \bigr) = 1 . $$

\ss

Basta mostrarmos agora que :
$$ \lim_{N \rightarrow \infty} \big\vert \bbP\bigl(R^{\eta}_N > t(N) \bigr) - \bbP\bigl( R_N > t(N) \bigr) \big\vert = 0 .$$


Mas pelas proposi\c c\~oes I e II temos que:

\begin{align*}
\big\vert \bbP\bigl( R^{\eta}_N > t(N) \bigr) &- \bbP \bigl( R_N > t(N) \bigr) \big\vert \cr 
&\leq \sup_{\eta} \bbP\bigl( \sigma^{\eta}(N^{1+\delta}) \not= \sigma^{+}(N^{1+\delta}) \bigr)  . \cr
\end{align*}

Pelo teorema I a probabilidade acima vai a 0 quando $ N $ diverge.

\hfill\newend

\noindent{\bf Obs.:} As proposi\c c\~oes III e IV nos fornecem limitantes para $ E(R_N) $ quando $ N $ diverge.

\ss

\ss

\noindent{\bf Lema III: } Considere $ \ \beta_N = \min( n \in \bbN : \bbP\bigl( R_N \geq n \bigr) \leq e^{-1} ) . $ Ent\~ao ,

$$\lim_{N \rightarrow \infty} \bbP\bigl( R_N \geq \beta_N \bigr) = e^{-1} . $$

\noindent{\bf Demonstra\c c\~ao:}

Pela defini\c c\~ao de $ \beta_N $ temos que:
$$ \bbP\bigl( R_N \geq \beta_N \bigr) \leq e^{-1} < \bbP\bigl( R_N \geq \beta_N -1
\bigr) $$

Como, 

$ 0 \leq \bbP\bigl(R_N \geq \beta_N -1 \bigr) - \bbP\bigl( R_N \geq \beta_N \bigr) \leq
\bbP\bigl( \beta_N -1 \leq R_N < \beta_N \bigr) $ ;

concluiremos  a demonstra\c c\~ao utilizando a propriedade de Markov.

\begin{align*}
\bbP\bigl( &\beta_N-1 \leq R_N < \beta_N \bigr) = \bbP\bigl( \beta_N - 1 \leq R_N < \beta_N \ | \ \sigma^{+}( \beta_N -1) \not\in F \bigr) \ \times \cr \cr
&\times \bbP\bigl( \sigma^{+}( \beta_N -1) \not\in F \bigr) \cr \cr
&= \bbP\bigl( \sigma(1) \in F \ | \ \sigma(0) \not\in F \bigr)  \ \bbP\bigl( \sigma^{+}( \beta_N -1) \not\in F \bigr) \cr 
&= \bbP\bigl( \sigma(1) \in F \ | \ \sigma(0) \not\in F \bigr) \ {\bf\Bigl[} \bbP\bigl(
\sigma^{+}( \beta_N -1) \not\in F \bigr) \ - \cr \cr
&- \sum_{\eta \in H_N} \nu(\eta) \bbP\bigl( \sigma^{\eta}( \beta_N -1) \not\in F \bigr) \ + \cr \cr
&+ \sum_{\eta \in H_N} \nu(\eta) \bbP \bigl( \sigma^{\eta}( \beta_N - 1) \not\in F \bigr) {\bf \Bigr]} \cr \cr
&= \bbP\bigl( \sigma(1) \in F \ | \ \sigma(0) \not\in F \bigr)  {\bf \Bigl[}
\bbP\bigl( \sigma^{+}( \beta_N -1) \not\in F \bigr) \ - \cr \cr
&- \sum_{\eta \in H_N} \nu(\eta) \bbP\bigl( \sigma^{\eta}( \beta_N -1) \not\in F {\bf \Bigr]} \ + \cr \cr
&+ { {2^{N} - |F|} \over {2^N}} \bbP\bigl( \sigma(1) \in F \ | \ \sigma(0) \not\in F \bigr)  \cr \cr
&\leq \bbP\bigl( \sigma(1) \in F \ | \ \sigma(0) \not\in F \bigr) \sum_{\eta \in H_N} \nu(\eta) {\bf \Bigl[} \bbP\bigl(\sigma^{+}(\beta_N -1) \not\in F \bigr) \ - \cr \cr
&- \bbP\bigl( \sigma^{\eta}(\beta_N -1) \not \in F \bigr) { \bf \Bigr]} + 
{ {2^{N} - |F|} \over { 2^{N}} } { {N^{\gamma}} \over N}  \cr \cr
&\leq \sum_{\eta \in H_N} \nu(\eta) \sup_{\eta} \bbP\bigl( \sigma^{+}(\beta_N-1) \not= \sigma^{\eta}(\beta_N-1) \bigr) + { {2^{N} - |F|} \over {2^{N} } } { {N^\gamma} \over {N}} \ . \cr
\end{align*}

Passando ao limite quando $N $ diverge temos que o primeiro termo vai a zero pelo teorema  I e o segundo vai a zero pelo corol\'ario III e pela hip\'otese.

Portanto,
$$ \lim_{N \rightarrow \infty} \bbP\bigl( R_N \geq \beta_N \bigr) = e^{-1} . $$

\hfill\newend

\noindent{\bf Nota\c c\~ao:} $ \sigma(\beta_N) = \sigma( [\beta_N] ) . $

\ss

\ss
\noindent{\bf Lema IV: } \ Existe um n\'umero real $ \alpha $ satisfazendo $ e^{-1} \leq \alpha< 1 $  tal que para $ N $ suficientemente grande e qualquer inteiro $ n $ temos:
$$ \bbP\bigl( R_N \geq \beta_N n \bigr) \leq \alpha^{n} . $$
\noindent{\bf Demonstra\c c\~ao:} \ A verifica\c c\~ao do resultado ser\'a feita por indu\c c\~ao.

\ss

Para $ n=1 $ o resultado \eh\ direto pois por defini\c c\~ao

$ \beta_N = \min\{n \in \bbN : \bbP\bigl( R_N \geq \beta_N \bigr) \leq e^{-1} \} . $

\ss

\ss

Assumiremos agora que a desigualdade vale para o inteiro $ n $. Provaremos para $ n+1 $ usando a propriedade de Markov.

\begin{align*}
\bbP\bigl( R_N \geq \beta_N (n+1) \bigr) &= \sum_{\eta \not\in F} \bbP\bigl(R_N \geq \beta_N n, \sigma(\beta_N n)=\eta \bigr) \times
\cr \cr
&\bbP\bigl(R_N \geq \beta_N \ | \ \sigma(0)=\eta \bigr) \cr \cr
&\leq \bbP\bigl(R_N \geq \beta_N n\bigr) \sup_{\eta \not\in F} \bbP\bigl( R_N \geq \beta_N \ | \ \sigma(0)=\eta \bigr) \cr \cr
&\leq \alpha^{n} \sup_{\eta \not\in F} \bbP\bigl(R_N \geq \beta_N \ | \
\sigma(0)=\eta \bigr) \cr
\end{align*}

Agora pela proposi\c c\~ao I e II temos que:

\begin{align*}
\big\vert \bbP\bigl( R^{\eta}_N \geq \beta_N \bigr) &- \bbP\bigl( R_N \geq \beta_N \bigr) \big\vert = \cr \cr
&\big\vert \bbP\bigl( R^{\eta}_N \geq \beta_N, \sigma^{\eta}(u) \not\in F, u\leq N^{1+\delta} \bigr) - \cr \cr
&- \bbP\bigl( R_N \geq \beta_N, \sigma^{+}(u) \not\in F, u\leq N^{1+\delta} \bigr) \big\vert \cr \cr
&\leq \sup_{\eta} \bbP\bigl( \sigma^{+}(N^{1+\delta}) \not= \sigma^{\eta}(N^{1+\delta}) \bigr) \cr
\end{align*}
Pelo teorema I temos que a probabilidade acima converge a zero quando $ N $ diverge.

Assim, para $ N $ suficientemente grande temos:

$$ \bbP\bigl(R_N > \beta_N(n+1) \bigr) \leq \alpha^{n} e^{-1} \leq \alpha^{n+1} . $$
\hfill\newend

\ss

\noindent{\bf Teorema III:} \ Para $ 0<\gamma<1 $  temos que:
$$ \lim_{N \rightarrow \infty} \bbP\Bigl( {{R_N } \over {\beta_N}} > t \Bigr) = e^{-t} . $$

\s
\noindent{\bf Demonstra\c c\~ao:} \ Para verificarmos que $ R_N $ normalizado por $ {\beta_N}  $ \ tem lei exponencial de par\^ametro 1 quando $ N $ diverge basta provarmos que :

\ss
\ss

\noindent{\bf a) } $ \ \lim\limits_{N \rightarrow \infty} { {E(R_N)} \over {\beta_N} } = 1. $

\ss

\ss

\noindent{\bf b) } $ \ \lim\limits_{N \rightarrow \infty} \big\vert \bbP\Bigl( R_N > { {\beta_N (t+s)}   }\Bigr) - \bbP\Bigl( R_N > {{\beta_N t} }\Bigr) \ \bbP\Bigl(
R_N > { {\beta_N s}  } \Bigr) \Big\vert = 0 , $ 

\ \ \ \ \ para qualquer $ s, t $ fixados.

\ss

\ss

\ss

\noindent{\bf obs.: 1)}  O segundo item garante que se a lei de $ { {R_N} \over { \beta_N}} $ converge quando $ N \rightarrow \infty $, o limite precisa ser uma lei exponencial ( talvez degenerada). \ Por outro lado, o lema III juntamente com este item implicam que se $ t $ \eh\ um n\'umero racional positivo, ent\~ao o limite 
$$ \lim_{N \rightarrow \infty} \bbP\bigl( R_N \geq \beta_N t \bigr) $$
existe e \eh\ igual a $ e^{-t} . $ Como a lei exponencial \eh\ cont\'\i nua, isto \eh\ suficiente para provar a converg\^encia para todo $ t $ real, o que conclui a prova.

\ss

\noindent{\bf \ \ \ \ \ 2)} Esta t\'ecnica foi utilizada em [2] para a demonstra\c c\~ao de outro resultado podendo servir como refer\^encia.

\ss

\ss

\noindent{\bf prova de a) : } 
\begin{align*}
{ {E(R_N)} \over {\beta_N} } &= { 1 \over {\beta_N} } \int_{0}^{\infty} \bbP\bigl( R_N>t \bigr) dt \cr \cr
&= \int_{0}^{\infty} \bbP\bigl( {{R_N} \over {\beta_N}} > t \bigr) dt .\cr
\end{align*}

\ss

\ss

\ss

Passando ao limite em $ N $ e utilizando o lema IV podemos aplicar o teorema da converg\^encia dominada de Lebesgue obtendo:

$$ \lim_{N \rightarrow \infty} \int_{0}^{\infty} \bbP\bigl( {{R_N} \over {\beta_N} } > t \bigr) dt = 1 . $$

\ss

\noindent{\bf prova de b) : } \ Primeiramente vamos mostrar o resultado para uma configura\c c\~ao inicial escolhida uniformemente em $ H_N $ . 

\ss

\noindent Observe os seguintes fatos:

\noindent{\bf Fato 1: }
\begin{align*}
 &\Big\vert \sum_{\eta \in H_N} \nu(\eta) \bbP\Bigl( R^{\eta}_{N}> \beta_N  (t+s) \Bigr) \cr \cr
 &- \sum_{\eta \in H_N} \nu(\eta) \bbP\Bigl( \sigma^{\eta}(u) \not\in F, \cr 
 & \ \ \ \ \ \ \ \ \ \ \ \ \forall u \in \{ 1,\dots,
 \beta_N t \} \cup \{ \beta_N t + N^{1 + \delta}, \dots, \beta_N (t+s) \} \Bigr) \Big\vert \cr \cr
 &\leq \sum_{\eta \in H_N} \nu(\eta) \bbP\Bigl( \sigma^{\eta}(u) \in F, \cr 
 & \ \ \ \ \ \ \ \ \ \ \ \hbox{ para algum } \ u \in
 \{ \beta_N t + 1,\dots,\beta_Nt + N^{1 + \delta} \} \Bigr) \cr
&\leq \sum_{\eta \in H_N} \nu(\eta) \sum_{u=\beta_N t}^{\beta_N t+N^{1+\delta}} \sum_{\zeta \in F} \bbP\Bigl( \sigma^{\eta}(u) = \zeta \Bigr) \cr \cr
&\leq { {N^{1+\delta} \bigl(N^{\gamma} +1 \bigr)} \over {2^N} } . \cr
\end{align*}
\ss
Assim a probabilidade acima vai a zero quando $ N $ diverge.

\ss

\noindent{\bf Fato 2: } 
\begin{align*}
&\Big\vert \sum_{\eta \in H_N} \nu(\eta) \bbP\Bigl( R^{\eta}_{N} >  \beta_N  s \Bigr) \cr \cr
&- \sum_{\eta \in H_N} \nu(\eta) \bbP\Bigl(\sigma^{\eta}(u) \not\in F, \forall u \in 
\{ N^{1 + \delta},\dots, \beta_N  s \} \Bigr) \Big\vert \cr \cr
&\leq \sum_{\eta \in H_N} \nu(\eta) \bbP\Bigl( \sigma^{\eta}(u) \in F, \hbox{para algum } u \in \{1,\dots,N^{1+\delta} \} \Bigr) \cr \cr
&\leq \sum_{\eta \in H_N} \nu(\eta) \sum_{u=1}^{N^{1+\delta}} \sum_{\zeta \in F} \bbP\Bigl( \sigma^{\eta} = \zeta  \Bigr)  \cr \cr
&\leq { {N^{1+\delta} (N^{\gamma}+1)} \over { 2^N} } . \cr
\end{align*}

Novamente, quando $ N $ diverge a probabilidade acima vai a zero.  

\ss

\noindent{\bf Fato 3: } \ Considere a seguinte express\~ao:
\begin{align*}
&\Big\vert \sum_{\eta \in H_N} \nu(\eta) \bbP\Bigl( R^{\eta}_{N} > \beta_N  (t+s) \Bigr)
\cr \cr
&- \sum_{\eta \in H_N} \nu(\eta) \bbP\Bigl( R^{\eta}_N > \beta_N  t \Bigr) \sum_{\eta \in H_N} 
\nu(\eta) \bbP\Bigl( R^{\eta}_{N} > \beta_N s \Bigr) \Big\vert . \cr
\end{align*}

Agora utilizando a propriedade de Markov, os fatos 1  e 2 a express\~ao acima \eh\ limitada por:
\begin{align*}
&\Big\vert \sum_{\eta \in H_N} \sum_{\kappa \not\in F} \nu(\eta)
\bbP\Bigl(R^{\eta}_{N}> \beta_N t, \sigma^{\eta}(\beta_N t) = \kappa \Bigr)
\cr \cr
&\times \bigl[ \bbP\Bigl(\sigma^{\kappa}(u) \not\in F, \forall u \in \{ N^{1+\delta},\dots,
\beta_N  s \} \Bigr) \cr \cr
&- \bbP\Bigl(\sigma^{\eta}(u) \not\in F, \forall u \in \{N^{1+\delta}, \dots, \beta_N  s \} \Bigr) \Bigr] \Big\vert \cr \cr
&\leq \sum_{\eta \in H_N} \nu(\eta) \sum_{\kappa \not\in F } \sup_{\kappa \in H_N}
\bbP\Bigl( \sigma^{\kappa}(N^{1 + \delta}) \not= \sigma^{\eta}(N^{1+\delta}) \Bigr) . \cr
\end{align*}

Assim, passando ao limite quando $ N $ diverge e utilizando o teorema I temos que a express\~ao acima vai a zero.

\ss

Resta-nos mostrar agora que:
$$\lim_{N \rightarrow \infty} \Big\vert \bbP\Bigl( R_{N} > \beta_N  t \Bigr) - 
\bbP\Bigl( R^{\eta}_{N} > \beta_N  t \Bigr) \Big\vert = 0 . $$

Com efeito,
\begin{align*}
\Big\vert \bbP\Bigl( R_{N}>\beta_N  t\Bigr)&-\bbP\Bigl( \sigma^{+}(u) \not\in F,
\forall u \in \{ N^{1+\delta},\dots,\beta_N t \} \Bigr) \Big\vert \cr \cr
&\leq \bbP\Bigl( R_{N}< N^{1+\delta} \Bigr) .\cr
\end{align*}

Da mesma forma,

\begin{align*}
 \Big\vert \bbP\Bigl(R^{\eta}_{N} > \beta_N t \Bigr)&-\bbP\Bigl(\sigma^{\eta}(u) \not\in F,
 \forall u \in \{ N^{1+\delta},\dots,\beta_N t \}\Bigr) \Big\vert \cr \cr
 &\leq \bbP\Bigl(R^{\eta}_{N}<N^{1+\delta} \Bigr) . \cr
\end{align*}

Portanto,
\begin{align*}
 \Big\vert \bbP\Bigl(R_{N} > \beta_N t\Bigr) &- \bbP\Bigl(R^{\eta}_{N}> \beta_N t \Bigr) \Big\vert \cr \cr
 &\leq \big\vert \bbP\Bigl( \sigma^{+}(u) \not\in F, \forall u \in \{N^{1+\delta},\dots,\beta_N t \} \Bigr) \cr \cr
 &- \bbP\Bigl( \sigma^{\eta}(u) \not\in F, \forall u \in
 \{ N^{1+\delta},\dots,\beta_N  t \} \Bigr) \big\vert \cr \cr
 &\leq \sup_{\eta \in H_N} \bbP\Bigl( \sigma^{+}(N^{1+\delta}) \not= \sigma^{\eta}(N^{1+\delta}) \Bigr) . \cr
\end{align*}

Agora passando ao limite quando $ N $ diverge e utilizando o teorema I temos o resultado.

\hfill\newend
\ss
\vfill\eject

\noindent{\bf Lema V : } ( Lema da Reflex\~ao )

\s

Para $ u , s $ fixados temos:
$$ \bbP \Bigl( V[0,s] \cap V[ s+1, s+u ] = \emptyset \Bigr) =
\bbP \Bigl( V[ 0,u-1 ] \cap V[ u,u+s ] = \emptyset \Bigr) .$$

\noindent{\bf Demonstra\c c\~ao : }

\ss

Seja $ {\cal C} $ a classe de conjuntos tal que:
$${\cal C } = \bigl\{ (i_1,\dots,i_{2c}); c=1,2,\dots : \sum\limits_{k=1}^{2c}
{\bf 1}_{\{ i = i_k \} } \ \hbox{ \eh\ par para todo } \ i \in \{ 1,\dots,N\} \bigr\} .$$
Note que se $ \sigma(k) = \eta^{i_1 \dots i_k} = \eta $ \ ent\~ao \ $ ( i_1, \dots, i_k ) \in {\cal C } $ .

\ss

Primeiramente notamos que:
\begin{align*}
\bbP\Bigl( V[0,s] &\cap V[s+1,s+u] = \emptyset\Bigr) \cr 
&= \bbP\Bigl( \eta \not\in V[s+1,s+u], \dots, \sigma(s) \not\in V[s+1,s+u] \Bigr) \cr 
&= \bbP\Bigl( (i_1,\dots,i_{m_1}) \not\in {\cal C}, \hbox{\ para } \ s+1 \leq m_1 \leq s+u, \dots , \cr 
&\dots \ \ \ (i_s,\dots,i_{m_s}) \not\in {\cal C}, \hbox{ \ para } \ s+1 \leq m_s \leq s+u \Bigr) \cr
\end{align*}

Por outro lado,
\begin{align*}
\bbP\Bigl( V[0,u-1] &\cap V[u,u+s] = \emptyset \Bigr) \cr 
&= \bbP\Bigl( \sigma(u) \not\in V[0,u-1],\dots, \sigma(u+s) \not\in V[0,u-1] \Bigr) \cr 
&= \bbP\Bigl( (i_{m_1},\dots,i_u) \not\in {\cal C }, \ \hbox{para \ } 1 \leq m_1 \leq u-1, \dots, \cr 
&\dots (i_{m_s},\dots,i_{u+s}) \not\in {\cal C} \ \hbox{para \ } 1 \leq m_s \leq u-1 \Bigr) \cr
\end{align*}

\s

Portanto temos a igualdade das probabilidades.  

\s

\hfill\newend

\ss

\noindent{\bf Corol\'ario V : } 

\s

$$ \lim_{N \rightarrow \infty} \bbP \Bigl(
V[ 0,N^{a} ] \ \cap \ V[ N^{a} + 1,N^{a} + N^{\gamma}+1 ] \not= \emptyset \Bigr) = 0 .$$

\ss

\noindent{\bf Demonstra\c c\~ao :}

\ss

Utilizando o lema IV temos que:

\s

\begin{align*}
&\bbP \Bigl( V[0,N^{a}] \cap V[ N^{a}+1,N^{a} + N^{\gamma} +1] = \emptyset \Bigr) 
\cr \cr
&= \bbP \Bigl( V[ 0,N^{\gamma} ] \cap V[ N^{\gamma}+ 1
, N^{a} + N^{\gamma} + 1 ] = \emptyset
\Bigr) \cr \cr
&= \bbP \Bigl( R_N > N^a + N^\gamma + 1 \Bigr) . \cr
\end{align*}

Passando ao limite e utilizando a proposi\c c\~ao IV temos:

\s

$$ \lim_{N \rightarrow \infty} \bbP\Bigl( V[0,N^a] \cap V[N^{a} + 1,N^{a} + N^{\gamma} + 1] \not= \emptyset \Bigr) = 0 . $$

\hfill\newend

\ss

\ss

\ss

\ss

\newpage
\def\eh{\rlap{e}\'{ }}
\def\Eh{\rlap{E}\'{ }}
\def\bbR{{I\kern-0.3emR}}
\def\bbP{{I\kern-0.3emP}}
\def\bbC{{I\kern-0.6emC}}
\def\bbN{{I\kern-0.3emN}}
\def\newend{\vbox to 7pt{\hbox to 7pt{\vrule height 7pt width 7pt}}}
\baselineskip = 18pt
\def\ss{\smallskip\smallskip}
\def\s{\smallskip}
\hyphenation {pro-ba-bi-li-da-de}
\hyphenation {a-le-a-to-ri-as}
\def\bff#1{\setbox0=\hbox{$#1$}
	\kern-.025em\copy0\kern-\wd0
	\kern.05em\copy0\kern-\wd0
	\kern-.025em\raise.0433em\box0 }

\noindent $ \bff{ \hbox{{\bf III.4} \underbar{\sl { Demonstra\c c\~ao do teorema IV } }}} $

\ss

\ss

\ss

Seja  $ \Theta = \inf \Bigl( t > 0; \xi(t) \in M \Bigr) ; $
ou seja, o tempo que o processo $ \xi(t) $ \ 
leva para alcan\c car o conjunto $ M . $

\ss
\ss

\noindent{\bf Proposi\c c\~ao V :}  \ \ Para $ 0   <  \gamma  <  1, $ \ temos:

\s
$$  \lim_{N \rightarrow \infty} \overline{E} \bigl( \ \bbP ( \Theta >
N^\gamma t \ ) \ \bigr) \ 
= \ e^{-t} . $$

\ss

\noindent{\bf Demonstra\c c\~ao : } \ \ \ Notamos primeiro que para
$ \overline{\omega} \in \overline{\Omega} $ fixado :

\ss

\ss

$ {\bbP} \Bigl( \Theta > {N^\gamma} t \Bigr) =
{\bf 1}_{ \{ \eta \not\in M \} }
{ 1 \over N } \sum\limits_{i_1 = 1}^N  {\bf 1}_{ \{ \eta^{i_1}
	\not\in M \} } \dots
{ 1 \over N } \sum\limits_{i_{{N^\gamma} t} = 1}^N
{\bf 1}_{ \{ \eta^{i_1 \dots i_{{N^\gamma} t}}  \not\in M \} } .$

\ss
\ss

\noindent Denotaremos por

\ss

\noindent {$ F_1(N) = \{ (i_1,\dots,i_{N^{\gamma}t}) \in \{ 1,\dots,N \}^{N^{\gamma}t} , 
	\ \forall \ l=1, \dots,N^{\gamma}t, \eta^{i_1 \dots i_l} \not\in \{ \eta,\eta^{i_1},\dots,\eta^{i_1 \dots i_{l-1}} \} \} ; $ \ \

	\noindent e por $ \ F_2(N) =
	\{1,\dots,N\}^{N^{\gamma}t} \backslash F_1(N).
	$  }

Seja $ i^{\prime} = (i_1, \dots,i_{N^{\gamma}t}) $ ent\~ao,

\begin{align*}
\overline{E} \Bigl( \bbP( \Theta > N^{\gamma}t ) \Bigr) &=
{ 1 \over {N^{ {N^\gamma}t } }} \sum_{ i^{\prime} \in F_1(N) }  \ \overline{E}( {\bf 1}_{ \{ \eta \not\in M \} } \dots
{\bf 1}_{ \{ \eta^{i_1 \dots i_{N^{\gamma}t}}  \not\in M \} })  \cr \cr
&+  { 1 \over { N^{{N^\gamma}t }}}  \sum_{ i^{\prime} \in F_2(N)} \ \overline{E} (  {\bf 1}_{ \{ \eta \not\in M \} } \dots
{\bf 1}_{ \{ \eta^{i_1 \dots i_{N^{\gamma}t}}  \not\in M \} }) \cr \cr
&=  {{|F_1(N)|} \over {N^{{N^\gamma}t}}}  \ ( 1 - { 1 \over { N^{\gamma}}})^{N^{\gamma}t + 1}  \cr \cr
&+  { 1 \over  {N^{{N^\gamma}t} }} \sum_{ i^{\prime} \in F_2(N)} \overline{E} (  {\bf 1}_{ \{ \eta \not\in M \} } \dots
{\bf 1}_{ \{ \eta^{i_1 \dots i_{N^{\gamma}t}}  \not\in M \} }) \cr
\end{align*}

\ss

\ss

\vfill\eject

Pelo corol\'ario III,

$$ \lim_{N \rightarrow \infty} { { |F_1(N)|} \over {N^{N^{\gamma}t}}} = 1 , \  \  \ \hbox{e} $$

\s

$$ \lim_{N \rightarrow \infty} {{ |F_2(N)|} \over {N^{N^{\gamma}t}}} = 0 . $$

Portanto,

$$\lim\limits_{N \rightarrow \infty} \overline{E} \Bigl( \bbP ( \Theta > 
N^{\gamma}t ) \Bigr) = \lim\limits_{N \rightarrow \infty} 
\Bigl( 1- { 1 \over {N^{\gamma}t}} \Bigr)^{N^{\gamma}t + 1} = e^{-t} .$$

\hfill \newend

\ss

\noindent{\bf Teorema IV :} \ \ \ Para $  \ 0  <  \gamma  <  1,
\  $ \  e $
\ \ \epsilon > 0 , $ \ temos :

$$ \lim_{N \rightarrow \infty } \overline{\bbP} \Bigl( \
\vert \bbP ( \Theta > { N^{\gamma} t} ) - e^{-t} \vert > \epsilon \Bigr) = 0 \ .$$

\ss

\noindent{\bf Demonstra\c c\~ao : }

\ss

Utilizando a desigualdade cl\'assica de Tchebyshev temos :

\begin{align*}
 \overline{\bbP} \Bigl( \vert {\bbP}
 ( \Theta > { N^\gamma } t ) - e^{-t} \vert  \ > \ \epsilon  \ \Bigr)
 &\leq  \ { 1 \over {\epsilon}^2} \ 
 \overline{E } \Bigl[ \Bigl( {\bbP} ( \Theta > { N^\gamma } t ) - e^{-t} \Bigr)^2
 \Bigr] \cr  \cr
 &= {1 \over {\epsilon}^2 } \bigl[ \ \overline{E } \Bigl[ \Bigl( {\bbP} ( \Theta >
 {N^\gamma} t ) \Bigr)^2 \Bigr]  \cr \cr
 &- 2 e^{-t} \overline{E} 
 \Bigl( \ \bbP ( \Theta > {N^\gamma} t ) \ \Bigr)
 + \overline{E} (e^{-t} )^2 \bigr] \ \cr \cr
 &= { 1 \over {\epsilon}^2 } 
 \bigl[ \ \overline{E} \Bigl[ \Bigl( \bbP ( \Theta > {N^\gamma} t ) \ \Bigr)^2
 \Bigr] \cr \cr
 &- 2 e^{-t} \overline{E} \Bigl(  \bbP ( \Theta > {N^\gamma} t ) \ \Bigr)
 +  e^{-2t}  \bigr] \ \cr
\end{align*}

\vfill\eject

Resta-nos calcular $ \overline{E} \Bigl[ \Bigl( {\bbP}
( \Theta > {N^\gamma} t ) \Bigr)^2 \Bigr] .$ 

\ss

\ss

Para  $ \overline{\omega} \in \overline{\Omega} $ fixado,

\begin{align*}
\bbP \bigl[\bigl( \Theta > N^{\gamma}t \bigr)^2\bigr] &= \Bigl\{
{\bf 1}_{ \{ \eta \not\in M \} } {1 \over N} \sum_{i_1=1}^{N} {\bf 1}_{ \{ \eta^{i_1} \not \in M \} } \dots {1 \over N } \sum_{i_{N^{\gamma}t}=1}^{N}
{\bf 1}_{ \{ \eta^{i_1 \dots i_{N^{\gamma}t}} \not\in M \} } \Bigr\} \times \cr \cr
&\times \Bigl\{{\bf 1}_{ \{ \eta \not\in M \} } { 1 \over N } 
\sum_{i^{*}_1 = 1}^{N} {\bf 1}_{ \{ \eta^{i^{*}_1} \not\in M \} } \dots
{1 \over N} \sum_{i^{*}_{N^{\gamma}t}=1}^{N} {\bf 1}_{ \{ \eta^{i^{*}_1 \dots i^{*}_{N^{\gamma}t} } \not\in M \} } \Bigr\} \cr
\end{align*}

\ss
\ss

\ss

Agora considere $ G^{*}= \bigl\{ \eta^{i^{*}_1}, \dots, \eta^{ i^{*}_1 \dots i^{*}_{N^{\gamma}t}} \bigr\} $ . Pelo corol\'ario III,
$$\lim_{N \rightarrow \infty} |G^{*}| - N^{\gamma}t = 0 .$$

Por outro lado, fixado $ \bigl( i^{*}_1,\dots,i^{*}_{N^{\gamma}t} \bigr) $ temos pelo teorema III que:

$$\lim_{N \rightarrow \infty} \bbP \Bigl( \{ \eta,\dots,\eta^{i_1 \dots i_{N^{\gamma}t}} \} \cap \{ \eta^{i^{*}_1}, \dots, \eta^{i^{*}_1 \dots i^{*}_{N^{\gamma}t}} \} \not= \emptyset \Bigr) = 0 . $$

Assim,

$$ \overline{E} \bigl[ \bbP(\Theta>N^{\gamma}t)^2 \bigr] = \bigl( 1 - 
{ 1 \over {N^\gamma}}\bigr)  \bigl( 1 - {1 \over {N^\gamma}}\bigr)^{2|G^{*}| } + o(N) . $$

Passando ao limite em $ \ N \ $ temos que :

$$ \lim_{N \rightarrow \infty } \overline{ E } \Bigl[
\Bigl( \bbP ( \Theta > N^\gamma t ) \Bigr)^2 \Bigr] = e^{-2t} . $$

Portanto, com o resultado da proposi\c c\~ao V :

$$ \lim_{N \rightarrow \infty} \overline{\bbP} \ \bigl[ \
\vert \bbP ( \Theta > N^\gamma t ) - e^{-t} \vert > \epsilon \ \bigr] \ 
= \ 0  . $$

\hfill \newend

\ss

\ss

\ss

\ss
\bibliographystyle{plain}
\nocite{tese1,tese2}
\bibliography{REF}

\end{document}